\documentclass[final,leqno]{siamltex}
\usepackage[latin1]{inputenc}
\usepackage[T1]{fontenc}
\usepackage[english]{babel}
\usepackage{times}
\usepackage{amssymb} 
\usepackage[pdftex]{graphicx} 
\usepackage{tikz}

\usepackage{theorem}
\usepackage{multicol}
\usepackage{float}
\usepackage{amsmath}
\usepackage{amsfonts}
\usepackage{hyperref}
\usepackage{setspace}
\usepackage{xcolor}

\usepackage[ruled,linesnumbered,resetcount]{algorithm2e}

\newcommand{\ii}{{\mathrm{i}}}
\newcommand{\sigmab}{\boldsymbol\sigma}
\newcommand{\mub}{\boldsymbol\mu}

\newcommand{\rhob}{\boldsymbol\rho}
\newcommand{\omegab}{\boldsymbol\omega}
\newcommand{\bm}[1]{\mathbf{#1}}
\newcommand{\R}{{\mathbb{R}}}
\newcommand{\C}{{\mathbb{C}}}
\DeclareMathOperator{\Imag}{Im}
\DeclareMathOperator{\Real}{Re}

\title{FDEMtools: a MATLAB package for FDEM data inversion} 

\author{G.P. Deidda, P. D\'iaz de Alba, C. Fenu, G. Lovicu and G. Rodriguez}

\begin{document}
\maketitle

\begin{abstract}
Electromagnetic induction surveys are among the most popular techniques for
non-destructive investigation of soil properties, in order to detect the
presence of both ground inhomogeneities and particular substances. This paper
introduces a MATLAB package, called FDEMtools, for the inversion of frequency
domain electromagnetic data collected by a ground conductivity meter, which
includes a graphical user interface to interactively modify the parameters of
the computation and visualize the results. Based on a nonlinear forward model
used to describe the interaction between an electromagnetic field and the soil,
the software reconstructs the distribution of either the electrical
conductivity or the magnetic permeability with respect to depth, by a
regularized damped Gauss--Newton method. The regularization part of the
algorithm is based on a low-rank approximation of the Jacobian of the nonlinear
model. The package allows the user to experiment with synthetic and
experimental data sets, and different regularization strategies, in order to
compare them and draw conclusions.
\end{abstract}



\pagestyle{myheadings}
\thispagestyle{plain}

\section{Introduction}
\label{sec:1}

Electromagnetic induction (EMI) techniques are often used for non-destructive
investigation of soil properties which are affected by electromagnetic
features of the subsurface layers, namely the electrical conductivity
$\sigma$ and the magnetic permeability $\mu$.
Knowing such parameters allows one to ascertain the presence of particular
substances, which is essential in many important applications: hydrological
characterizations~\cite{bgddrc18,Cassiani}, hazardous waste
studies~\cite{Martinelli}, archaeological surveys~\cite{Lascano,Osella},
precision agriculture~\cite{dcadlrvc18,Yao}, unexploded ordnance
detection~\cite{Huang}, etc.

A ground conductivity meter is the typical measuring instrument for
frequency domain electromagnetic (FDEM) induction techniques.
It is composed by two coils (a transmitter and a receiver)
placed at a fixed distance. The dipoles may be aligned vertically or
horizontally with respect to the ground level. An alternating sinusoidal
current in the transmitter produces a primary magnetic field $H_P$, which
induces small eddy currents in the subsurface. These currents produce, in turn,
a secondary magnetic field $H_S$, which is sensed by the receiver. The 
ratio $H_S/H_P$ of the secondary to the primary magnetic fields is measured by
the device, providing information about the amplitude and the
phase of the signal. The real part, or the in-phase component, of the measured
signal is mainly affected by the magnetic permeability of the subsoil; the
imaginary part, also called the out-of-phase or quadrature component, mainly by
the electrical
conductivity.

In this work, we present a MATLAB package for the numerical inversion of a
nonlinear model which describes the interaction between an electromagnetic
field and the soil. The computation is performed by a Gauss--Newton method and
regularized by means of either the truncated singular value decomposition
(TSVD), or the generalized truncated singular value decomposition (TGSVD).
The computation of the forward model, as well as the analytical expression of
its Jacobian matrix, is performed by some functions in the package. Either the
electrical conductivity or the magnetic permeability of the soil with respect
to depth can be reconstructed, if an assumption can be made on the behavior of
the other quantity.

The software development started during the research work which lead
to~\cite{dfr14}. In this paper, assuming the permeability of the soil to be 
known, the aim was to reconstruct the electrical conductivity taking as
input only the quadrature component of the measurements.
Here, the loop-loop device used a single scanning frequency, and multiple
measurements were obtained by varying the height of the instrument above the
ground and the orientation of the two coils. An inversion method was proposed,
based on a low-rank approximation of the Jacobian of the forward model.
Also, the analytical expression of the Jacobian with respect to the electrical
conductivity was computed for the first time, and it was compared to its finite
difference approximation for what regards accuracy and computing time.

The algorithm was extended to deal with multiple-frequency data sets
in~\cite{dr16}. This approach was motivated by the availability of devices
which can take simultaneous readings at multiple scanning frequencies. 
In the same paper, experiments were performed taking as input either the
in-phase or the quadrature components of the signal, in order to investigate
the possibility of extracting further information from the available data.

Later, the paper~\cite{ddr17} focused on the identification of the magnetic
permeability. The main result was to obtain stable analytical formulas for the
computation of the Jacobian of the forward model with respect to the variation
of the magnetic permeability.
The paper numerically investigated the conditioning of the problem, as well as
the performance of the analytical expression of the Jacobian in comparison to
its finite difference approximation.
Then, under the working assumption that the conductivity was known in advance,
the reliability of the inversion algorithm was verified, and the information
content of data sets obtained by varying either the scanning frequency of the
device or its height above the ground was investigated.
The generated data sets, contaminated by a noise level compatible with
real-world applications, allowed the authors to analyze the behavior of the
algorithm in a controlled setting.

Suitably extended version of the software were used in~\cite{ddrv18} to
investigate the effect of a particular sparsity enhancing regularization
technique, and in~\cite{bgddrc18,dcadlrvc18} to process specific experimental
data sets.
The algorithm was finally adapted in \cite{ddrv19} to the inversion of the
whole complex signal, rather than just one of its components. The resulting
software, that is essentially the one presented in this paper, was applied to
an experimental data set collected at the Molentargius Saline Regional Park, Sardinia, Italy \cite{ddrv19}.

The above mentioned references focus on a nonlinear model described, e.g., in
\cite{Ward}, and further studied in \cite{Hendrickx}. We remark that when the
soil conductivity takes relatively small values (below 0.5 S/m) a linear model
can be used to predict the measurements produced by a particular device
\cite{McNeill}.
A numerical method for this model was initially proposed in \cite{Borchers}.
The same model has recently been studied from the theoretical point of view in
\cite{dfrv19}, where an optimized solution method was proposed too.

The present paper is structured as follows. In Section \ref{sec:2} we give a
brief description of the nonlinear forward model and we introduce the inversion
algorithm, together with the regularization procedure. Section \ref{sec:3}
presents the FDEMtools package as well as its graphical user interface (GUI),
describing the installation process and how to use the software by means of
some numerical examples.
For extensive numerical experiments, we refer to the above mentioned papers.
The final Section \ref{sec:4} summarizes the content of the paper.

\section{Computational methods}\label{sec:2}

\subsection{The nonlinear forward model}
\label{subsec:1}

A forward model which predicts the data measured by a FDEM induction device,
when the distribution of the conductivity and the magnetic permeability in the
subsoil is known, has been described in~\cite{hendr02}. 


The model assumes the soil to be layered, so that both the electrical
conductivity and the magnetic permeability are piecewise constant functions. 
Each subsoil layer, of thickness $d_k$ (m), is characterized by an electrical
conductivity $\sigma_k$ (S/m) and a magnetic permeability $\mu_k$ (H/m), for
$k=1,\ldots,n$; see \cite{ddr17}. 
The thickness of the deepest layer $d_n$ is considered infinite.
The two coils of the measuring device are at height $h$ above the ground, their
distance is $\rho$.

Let $u_k(\lambda) = \sqrt{\lambda^2 + \ii\sigma_k\mu_k\omega}$, where $\omega$
is the angular frequency of the instrument, that is, $2\pi$ times the frequency
in Hertz. The variable $\lambda$ ranges from zero to infinity, 
it measures the ratio between the depth below the ground surface and the
inter-coil distance $\rho$.
If we denote the characteristic admittance in the $k$-th layer by $N_k(\lambda)
= u_k(\lambda)/(\ii\mu_k\omega)$, then it is shown in~\cite{wait2} that the
surface admittance $Y_k(\lambda)$ at the top of the same layer verifies the
recursion 
\begin{equation}
\label{surfadm}
Y_k(\lambda) = N_k(\lambda)\frac{Y_{k+1}(\lambda)+N_k(\lambda)
\tanh(d_k u_k(\lambda))}{N_k(\lambda) + Y_{k+1}(\lambda)
\tanh(d_k u_k(\lambda))},
\end{equation}
for $k=n-1,\ldots,1$. 
For $k=n$, the characteristic admittance and the surface admittance are assumed
to coincide, that is, the value $Y_n(\lambda) = N_n(\lambda)$ is used to
initialize the recursion.
Both the characteristic and the surface admittances are functions
of the frequency $\omega$ via the functions $u_k(\lambda)$.

Let us define the \emph{reflection factor} as
\begin{equation}\label{reflfact}
R_{\omega,0}(\lambda) = \frac{N_0(\lambda) - Y_1(\lambda)}{N_0(\lambda) + Y_1(\lambda)},
\end{equation}
where $Y_1(\lambda)$ is computed by the recursion \eqref{surfadm}.
The ratio of the secondary to the primary field for the vertical and horizontal
orientation of the coils, respectively, are given by
\begin{equation}
\begin{aligned}
M^V(\sigmab,\mub;h,\omega,\rho) &= -\rho^3 \mathcal{H}_0\left[
\lambda e^{-2h\lambda} R_{\omega,0}(\lambda) \right](\rho), \\
M^H(\sigmab,\mub;h,\omega,\rho) &= -\rho^2 \mathcal{H}_1\left[
e^{-2h\lambda} R_{\omega,0}(\lambda) \right](\rho),
\end{aligned}
\label{nonlinmodel}
\end{equation}
where $\sigmab=(\sigma_1,\ldots,\sigma_n)^T$, $\mub=(\mu_1,\ldots,\mu_n)^T$,
$N_0(\lambda)=\lambda/(\ii\mu_0\omega)$, $\mu_0=4\pi\cdot 10^{-7}H/m$ is the
magnetic permeability of free space, and $R_{\omega,0}(\lambda)$ is defined by
\eqref{reflfact}.
We denote by
$$
\mathcal{H}_\nu[f](\rho) = \int_{0}^{\infty} f(\lambda) J_\nu(\rho\lambda)
\lambda \,d\lambda , \qquad \nu=0,1,
$$
the Hankel transform, where $J_0,J_1$ are first kind Bessel functions of order
0 and 1, respectively.

We remark here that the functions in \eqref{nonlinmodel} take complex values
and that the measuring devices, in general, return both the real
(\emph{in-phase}\/) and the imaginary (\emph{quadrature}\/) components of the
fields ratio.
The quadrature component, suitably scaled, can be interpreted as an
\emph{apparent conductivity}, while the in-phase component is related to the
magnetic permeability of the soil.

\subsection{Inversion procedure}
\label{subsec:2}

Simultaneous measurements with different inter-coil distances or different
operating frequencies can be recorded by recent FDEM devices at different
heights.
We denote by $\rhob=(\rho_1,\ldots,\rho_{m_\rho})^T$,
$\bm{h}=(h_1,\ldots,h_{m_h})^T$, and
$\omegab=(\omega_1,\ldots,\omega_{m_\omega})^T$, the vectors containing the
loop-loop distances, the heights, and the angular frequencies at which the
readings were taken.
We consider the corresponding $m=2 m_\rho m_h m_\omega$ data points
$b^\nu_{tij}$, where $t=1,\ldots,m_\rho$, $i=1,\ldots,m_h$,
$j=1,\ldots,m_\omega$, while $\nu=0,1$ represents the vertical and
horizontal orientations of the coils, respectively.
The observations $b^\nu_{tij}$ are rearranged in a vector $\bm{b}\in\C^m$.

Let us consider the complex residual vector 
$$
\bm{r}(\sigmab,\mub) = \bm{b}-\bm{M}(\sigmab,\mub;\bm{h},\omegab,\rhob)
$$
as a function of the conductivities $\sigma_k$ and the permeabilities $\mu_k$,
$k=1,\ldots,n$; the vector function $\bm{M}$ returns the readings predicted by
the model \eqref{nonlinmodel} in the same order the values $b^\nu_{tij}$
were arranged in the vector $\bm{b}$.
To process the data set, we solve the following nonlinear least-squares problem
\begin{equation}\label{minimization}
\min_{\sigmab,\mub\in\R^n} \frac{1}{2} \|\bm{r}(\sigmab,\mub)\|^2,
\end{equation}
by the Gauss--Newton method.
At each step of the iterative algorithm we minimize
the 2-norm of a linear approximation of the residual
$\bm{r}(\sigmab_{k+1},\mub_{k+1})=\bm{r}(\sigmab_k+\bm{s}_k,\mub_k+\bm{t}_k)$,
namely,
\begin{equation}\label{gaussnewt}
\displaystyle \min_\bm{q} \|\bm{r}(\sigmab_{k}, \mub_k) + J(\sigmab_{k}, \mub_k)\bm{q}_k\|,
\end{equation}
where 
\[
\bm{q}_k=\begin{bmatrix} \bm{s}_k \\ \bm{t}_k \end{bmatrix}, \qquad
\bm{s}_k,\bm{t}_k\in\C^n,
\]
and 
\[
J(\sigmab_k,\mub_k)= \begin{bmatrix} J^\sigma, & J^\mu\end{bmatrix}, \qquad
(J^\sigma)_{i,k} = \frac{\partial r_i(\sigmab, \mub)}{\partial \sigma_k},
\qquad 
(J^\mu)_{i,k} = \frac{\partial r_i(\sigmab, \mub)}{\partial \mu_k},
\]
for $i=1,\ldots,m$ and $k=1,\ldots,n$.

The residual function $\bm{r}$ being complex-valued, we solve 
problem \eqref{gaussnewt} by stacking the real and imaginary part of the
residual as follows \cite{ddrv19}
\begin{equation*}
\widetilde{\bm{r}}(\sigmab,\mub)= \left[\begin{array}{c}
\beta \cdot \Real(\bm{r}(\sigmab,\mub)) \\
\Imag(\bm{r}(\sigmab,\mub))
\end{array}\right] \in \R^{2m}, \quad
\widetilde{J}(\sigmab,\mub)=\left[\begin{array}{c}
\beta \cdot \Real(J(\sigmab,\mub)) \\
\Imag(J(\sigmab,\mub))
\end{array}\right] \in \R^{2m\times 2n}.
\end{equation*}
The positive parameter $\beta$ allows for balancing the contributions of
the in-phase and the quadrature component. It is in general set to one, unless
it is modified by the user on the basis of available a priori information,
e.g., on the soil composition or on the noise level.

So, we replace \eqref{gaussnewt} by
\begin{equation}
\min_{\bm{q} \in \R^{2n}}
\|\widetilde{\bm{r}}(\sigmab_k,\mub_k)+\widetilde{J}_k \bm{q}\|,
\label{eq8:leastsquaresnew}
\end{equation}
with $\widetilde{J}_k=\widetilde{J}(\sigmab_k,\mub_k)$,
and the iterative method becomes
$$
(\sigmab_{k+1},\mub_{k+1}) = (\sigmab_k + \alpha_k\bm{s}_k,\mub_k + \alpha_k\bm{t}_k) 
= (\sigmab_k,\mub_k) - \alpha_k \widetilde{J}_k^{\dagger} \, \widetilde{\bm{r}}(\sigmab_k,\mub_k),
$$
where $\widetilde{J}_k^{\dagger}$ is the Moore--Penrose pseudoinverse of $\widetilde{J}_k$ and
$\alpha_k$ is a damping parameter which ensures the convergence.
It is determined by coupling the Armijo--Goldstein principle~\cite{bjo96} to
the positivity constraint $(\sigmab_{k+1},\mub_{k+1})>0$ (componentwise)
\cite{ddr17,dfr14}.

The analytical expression of the Jacobian matrices with respect to the
electrical conductivity and the magnetic permeability were computed
in~\cite{dfr14} and~\cite{ddr17}, respectively. In fact, the same papers show
that the analytical Jacobian is more accurate and faster to compute, than
resorting to a finite difference approximation.

\subsubsection{Regularization}
\label{subsec:3}

It is well known that the minimization problem~\eqref{minimization} is
extremely ill-conditioned. In particular, it has been observed
in~\cite{ddr17,dfr14} that the Jacobian matrix $J$ has a large condition number
virtually for all values of the variables $\sigmab$ and $\mub$ in the solution
domain. A common remedy to address ill-conditioning consists of resorting to
regularization, i.e., replacing the linearized least-squares
problem~\eqref{eq8:leastsquaresnew} by a nearby problem, whose solution is less
sensitive to the propagation of the errors in the data.

A regularization method which particularly suits the problem, given the size of
the matrices involved, is the truncated singular value decomposition (TSVD).
The best rank $\ell$ approximation ($\ell \leq \kappa=\rank(\widetilde{J}_k)$) to the Jacobian
matrix, according to the Euclidean norm, is easily obtained by the SVD
decomposition. This factorization allows us to replace the ill-conditioned
Jacobian $\widetilde{J}_k$ by a well-conditioned low-rank matrix;
see~\cite{Hansen}.

When some kind of \emph{a priori} information for the problem is available,
e.g., the solution is a smooth function, it is useful to introduce
a regularization matrix $L\in\R^{p\times 2n}$ ($p\leq 2n$), whose null
space approximately contains the sought solution~\cite{rr13}. Under the
assumption $\mathcal{N}(\widetilde{J}_k)\cap\mathcal{N}(L)=\{\bm{0}\}$,
problem~\eqref{eq8:leastsquaresnew} is replaced by
\begin{equation}\label{min}
\min_{\bm{q}\in\mathcal{S}} \|L\bm{q}\|^2, \qquad
\mathcal{S} = \{ \bm{q}\in\R^{2n} \,:\, \widetilde{J}_k^T\widetilde{J}_k\bm{q}=-\widetilde{J}_k^T\widetilde{\bm{r}}(\sigmab_k, \mub_k) \}.
\end{equation}
Very common choices for $L$ are the discretization of the first or
second derivative operators.

Let the generalized singular value decomposition (GSVD)~\cite{gvl96} of the
matrix pair $(\widetilde{J}_k,L)$ be
$$
\widetilde{J}_k = U \Sigma_J Z^{-1}, \qquad L = V \Sigma_L Z^{-1},
$$
where $U$ and $V$ are matrices with orthogonal columns $\bm{u}_i$ and
$\bm{v}_i$, respectively, $Z$ is a nonsingular matrix with columns $\bm{z}_i$,
and $\Sigma_J$, $\Sigma_L$ are diagonal matrices with diagonal entries $c_i$
and $s_i$.
Under the assumption that $m=2 m_\rho m_h m_\omega<2n$, common in the
generality of cases, the truncated GSVD (TGSVD) solution $\bm{q}^{(\ell)}$
(see~\cite{Hansen} for details) can be written as
$$
\bm{q}^{(\ell)} = -\sum_{i=p-\ell+1}^p
    \frac{\bm{u}_i^T\bm{r}_k}{c_{i-2n+\kappa}}\, \bm{z}_i
    - \sum_{i=p+1}^{2n} (\bm{u}_i^T\bm{r}_k)\, \bm{z}_i,
$$
where $\kappa=\rank(\widetilde{J}_k)$, $\ell=1,\ldots,\kappa+p-2n$
is the regularization parameter, and $\bm{r}_k=\bm{r}(\sigmab_k,\mub_k)$.

The resulting regularized damped Gauss--Newton method reads
$$
(\sigmab_{k+1}^{(\ell)},\mub_{k+1}^{(\ell)}) =
(\sigmab_{k}^{(\ell)},\mub_{k}^{(\ell)}) + \alpha_k \bm{q}_k^{(\ell)},
$$
with $\ell$ fixed and $\alpha_k$ determined at each step as explained in
Section~\ref{subsec:2}.
We iterate until
$$
\|(\sigmab_k^{(\ell)},\mub_k^{(\ell)})-(\sigmab_{k-1}^{(\ell)},\mub_{k-1}^{(\ell)})\| < \tau \|(\sigmab_k^{(\ell)},\mub_k^{(\ell)})\|
\quad \text{or} \quad k > K_{\text{max}} \quad \text{or} \quad 
\alpha_k < \varepsilon.
$$
The solution at convergence is denoted by $(\sigmab^{(\ell)},\mub^{(\ell)})$.

The choice of the regularization parameter $\ell$ is crucial. In real
applications, experimental data are affected by noise. We express the data
vector in the residual function~\eqref{minimization} by
$\bm{b}=\widehat{\bm{b}}+\bm{e}$, where $\widehat{\bm{b}}$ contains the exact
data and $\bm{e}$ is the noise vector.
If the noise is Gaussian and an accurate estimate of $\|\bm{e}\|$ is available,
we determine $\ell$ by the discrepancy principle~\cite{Hansen,Morozov}.
When such an estimate is not at hand, we adopt heuristic methods, such as 
the L-curve, in the implementation presented in \cite{hjr07}, or an hybrid
method based on the quasi-optimality criterion, described in
\cite{morigi2006,rr13};
see~\cite{hrr15,prry18,rr13} for a review of similar methods.

Of course, the choice of the regularization matrix $L$ has strong effects on
the results, as it incorporates the available a priori information regarding
the solution. 
Consistently, whenever we have information about a \emph{blocky} behavior of
the solution, in order to maximize the spatial resolution of the result, we
consider an alternative regularization term called the \emph{minimum
gradient support} (MGS), which has been successfully applied in various
geophysical settings \cite{vfcka15,zhdanov02,zvu06}. This approach consists of
substituting the term $\|L\bm{q}\|^2$ in \eqref{min} by the nonlinear
stabilizing term \cite{ddrv19}
\begin{equation}
S_\tau(\bm{q}) = \sum_{r=1}^p \frac{\left( \frac{(L\bm{q})_r}{q_r}
\right)^2}{\left( \frac{(L\bm{q})_r}{q_r} \right)^2+\tau^2},
\label{sharp}
\end{equation}
where $L$ is a regularization matrix.

The nonlinear regularization term $S_\tau(\bm{q})$ favors the
sparsity of the solution and the reconstruction of blocky features.
Indeed, it can be shown \cite{pz99,vdc12} that, when $\tau$ approaches 0,
this approach minimizes the number of components where the
vector $L\bm{q}$ is nonzero.
Therefore, if $L$ is chosen to be the discretization of the first derivative
$D_1$, the stabilizer \eqref{sharp} selects the solution update corresponding
to minimal nonvanishing spatial variation; see~\cite{ddrv19} for further
comments and for a description of our implementation.

\section{The FDEMtools software package}
\label{sec:3}

The package consists of a set of MATLAB routines which implements the
algorithms for the inversion of FDEM data sketched in the previous sections. In
addition to the analysis and solution routines, the package includes 1D and 2D
test problems, showing how to process both synthetic and experimental data sets, and
investigate different device configurations and numerical strategies, in order
to compare them and draw conclusions.

The software is available from Netlib (http://www.netlib.org/numeralgo/) as the
na53 package. It requires P.~C.~Hansen's Regularization Tools package 
\cite{Hansentools} to be installed on the computer, and its directory to
be added to MATLAB search path.
FDEMtools is distributed as an archive file; by uncompressing it, a new
directory will be created, containing the functions of the toolbox as well as a
user manual. This directory must be added to MATLAB search path in order to be
able to use the software from other directories. More information on the
installation procedure can be found in the README.txt file, contained in the
toolbox directory.
All the routines are documented via the usual MATLAB \texttt{help} function.

Table~\ref{tab:funs} groups the MATLAB routines by subject area, giving a
short description of their purpose.
The first group ``Forward Model Routines'' includes the functions for computing
the forward model, that is, the model prediction for a given conductivity and
permeability distribution. The user can optionally compute also the Jacobian
matrix, by its analytical expression. The section ``Computational Routines''
describes the codes for the inversion algorithm, and the last section
``Auxiliary Routines'' lists some routines needed to complete the
whole process, that are unlikely to be called directly by the user.
See the file \texttt{Contents.m} for details.

\begin{table}[ht!]
\centering
\footnotesize \caption{FDEMtools reference.}
\begin{tabular}{p{2.4cm}p{8.5cm}}
\hline
\multicolumn{2}{c}{\textbf{Forward Model Routines}} \\
\hline
\texttt{reflfact} & compute the reflection factor $R_{\omega,0}(\lambda)$
	\eqref{reflfact} \\
\texttt{hratio} & compute the ratio $H_S/H_P$ \eqref{nonlinmodel}, i.e., the
	device readings \\
\texttt{inphase} & compute the in-phase (real) component of the ratio 
	$H_S/H_P$ \\
\texttt{quadracomp} & compute the quadrature (imaginary) component of 
	$H_S/H_P$ \\
\texttt{aconduct} & compute the apparent conductivity \\
\hline
\multicolumn{2}{c}{\textbf{Computational Routines}} \\
\hline
\texttt{emsolvenlsig} & reconstruct the electrical conductivity \\
\texttt{emsolvenlmu} & reconstruct the magnetic permeability \\
\texttt{tsvdnewt} & Gauss--Newton method regularized by T(G)SVD \\
\texttt{jack} & approximate the Jacobian matrix by finite differences \\
\texttt{hankelpts} & quadrature nodes for Hankel transform;
	see~\cite{anderson1979} \\
\texttt{hankelwts} & quadrature weights for Hankel transform;
	see~\cite{anderson1979} \\
\texttt{FDEMgui} & graphical user interface (GUI) activation command \\
\hline
\multicolumn{2}{c}{\textbf{Test Scripts}} \\
\hline
\texttt{driver} & test program for \texttt{emsolvenlsig} and \texttt{emsolvenlmu} \\
\texttt{driver2D} & 2D test program \\
\hline
\multicolumn{2}{c}{\textbf{Auxiliary Routines}} \\
\hline
\texttt{fdemcomp} & execute the inversion algorithm \\
\texttt{chooseparam} & define default parameters and test functions \\
\texttt{chooseparambis} & define default parameters \\
\texttt{fdemdoi} & compute the depth of investigation (DOI); see \cite{ddrv19} \\
\texttt{fdemprint} & print information about the whole process \\
\texttt{mgsreg} & compute the MGS regularization term; see \cite{ddrv19} \\
\texttt{addnoise} & add noise to data \\
\texttt{morozov} & choose regularization parameter by discrepancy principle \\
\texttt{quasihybrid} & choose regularization parameter by quasi-hybrid method; see \cite{rr13} \\
\texttt{plotresults} & display intermediate results during iteration \\
\texttt{fdemplot} & plot the reconstructed solution and, if available, the
	exact one \\
\end{tabular}
\label{tab:funs}
\end{table}

For example, the command
\begin{quote}
\footnotesize
\begin{verbatim}
[M,J] = hratio(sigma,mu,h,d,R,omega,'vertical','sigma');
\end{verbatim}
\end{quote}
computes the function $M^V$ in \eqref{nonlinmodel} (corresponding to the
``vertical'' orientation of the device), as well as its Jacobian with respect
to the $\sigma$ variables, under the assumption that the conductivity and
the permeability in each layer are defined, respectively, by the vectors
\texttt{sigma} and \texttt{mu}, and the layers thickness by the vector
\texttt{d}.
The variables \texttt{h}, \texttt{R}, and \texttt{omega} define the height of
the device above the ground, the inter-coil distance $\rho$, and the angular
frequency $\omega$, respectively; the last three variables may be vectors, in case of
multiple measurements.
The computed predictions are returned in a vector, according to a particular
ordering which is described in the documentation; see \texttt{help hratio}.
The reflection factor \eqref{reflfact} is computed by the function
\texttt{reflfact}, while the computation of the Hankel transform in
\eqref{nonlinmodel} is performed by a quadrature formula described in
\cite{anderson1979}, which is constructed by two routines from \cite{hendr02}.
The ``Forward Model'' functions are totally independent of the rest of the
package, so they could be used to test the performance of other inversion
algorithms.

The computational routines which are intended to be called directly by the user
are \texttt{emsolvenlsig} and \texttt{emsolvenlmu}.
The first one reconstructs the conductivity distribution with respect to depth,
taking as input the device measurements, and assuming the permeability
distribution is known. The second function does the same for the magnetic
permeability. 

The following example shows how to call the inversion function in order to
perform data inversion 
\begin{quote}
\footnotesize
\begin{verbatim}
% initialize mu,h,d,R,omega
func = @(sigma) hratio(sigma,mu,h,d,R,omega,'both','sigma');
opts = struct('func',func,'Lind',2,'orientation','both','show',1);
sigma = emsolvenlsig(func,y,h,d,opts);
\end{verbatim}
\end{quote}
After initializing some variables, the model function is chosen. In this case,
it is the complex signal in its entirety considering both orientations for the
device, i.e., both components of \eqref{nonlinmodel}. The permeability
\texttt{mu} is assumed to be known and the conductivity \texttt{sigma} is
determined by the regularized algorithm, assuming the data set is contained in
the vector \texttt{y}. Some options are set before the call and encoded in the
\texttt{opts} variable: the forward model, the regularization matrix (in this
case the discrete approximation of the second derivative), the orientation of
the coils (both vertical and horizontal), and the request to display
intermediate results.

In the case of multiple measurements, taken at different locations, each data
vector should be stored in a column of \texttt{y}. In this case the inversion
routine processes a column at a time, and returns a matrix containing in each
column the computed distribution of the electrical conductivity.

The package provides the \texttt{driver} program as an example of how to
perform the actual computations.
The first part of the \texttt{driver} defines the parameters which control the
construction of the synthetic data set, the device setting, the soil
discretization, the chosen model, and the tuning of the algorithm.
Then, the synthetic data set is constructed, it is contaminated by Gaussian
noise, and the inversion routine is called.

\begin{figure}[htbp]
\centerline{\includegraphics[width=.7\textwidth]{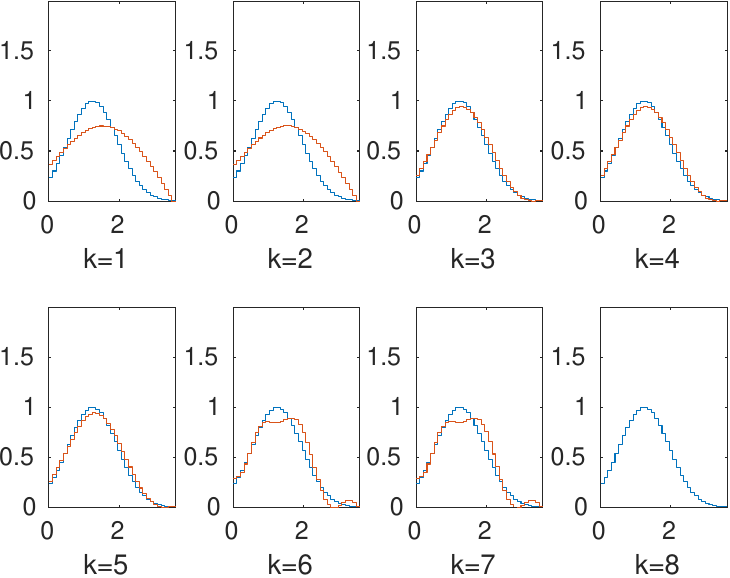}}
\caption{Regularized solutions computed by a call to \texttt{driver},
compared to the exact solution which generated the synthetic data set. The
regularization parameter, that is the TGSVD truncation parameter, is denoted by
$k$.}
\label{fig:iters}
\end{figure}

\begin{figure}[htbp]
\centerline{\includegraphics[width=.47\textwidth]{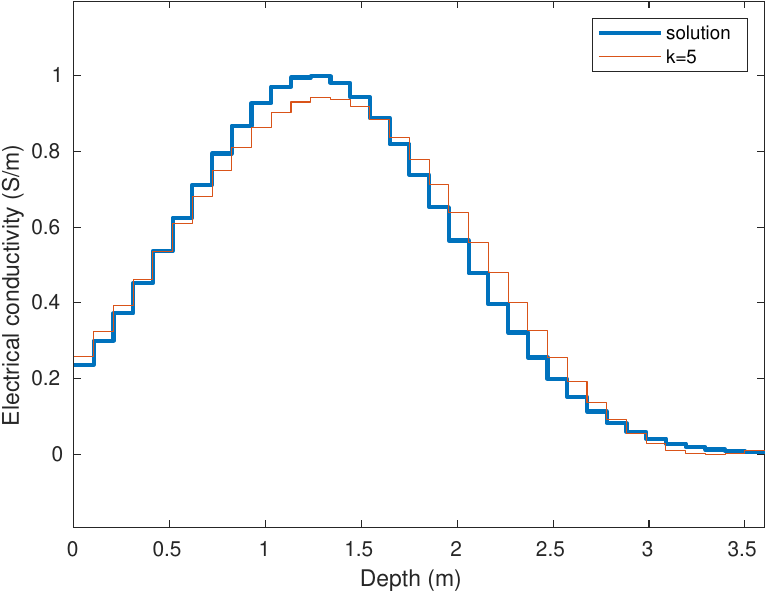}\hfill
\includegraphics[width=.5\textwidth]{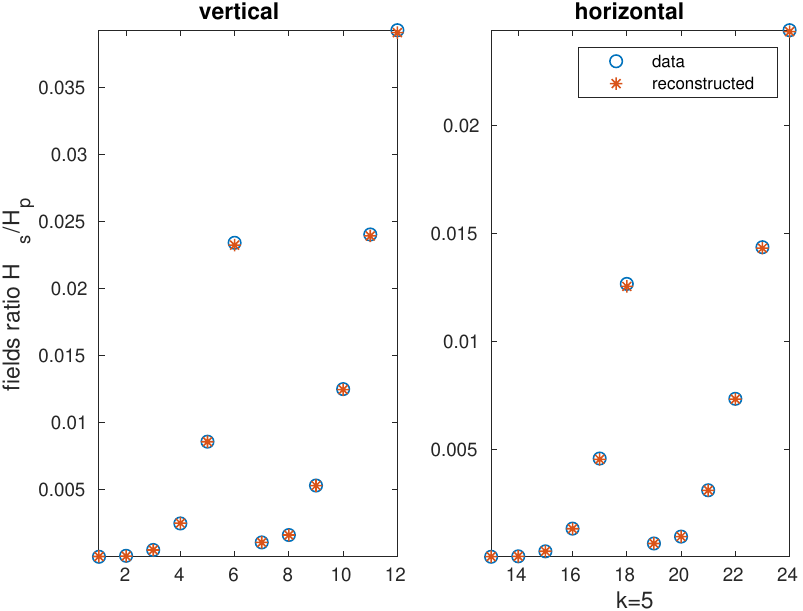}}
\caption{The graph on the left shows the solution selected by the ``corner''
criterion; the one on the right the data fitting between the synthetic data set
and the model prediction. This two graphs can be updated by clicking with the
left mouse button on one of the solutions displayed in
Figure~\ref{fig:iters}.}
\label{fig:soldata}
\end{figure}

The \texttt{driver} program provided with the package generates a synthetic
data set corresponding to a device operating with a single pair of coils, at 6
operating frequencies, and a single height above the ground. Both the vertical
and the horizontal orientations are considered so that there are $m=12$ data
points, but since the complex signal is processed the number of data values to
fit is doubled.
The subsoil is discretized in $n=35$ equal layers, up to a depth of 3.5m. We
assume a constant magnetic permeability $\mu_0$ (the value in the empty space)
and reconstruct the distribution of the electrical conductivity.

During the computation, information about the data set and the intermediate
results are displayed in the MATLAB main window and in some figures, which we
reproduce here.
Figure~\ref{fig:iters} shows the regularized solutions computed, while the
graph on the left of Figure~\ref{fig:soldata} displays the solution selected by
the ``corner'' criterion, and the one on the right the fitting between the
input data set and the model prediction corresponding to the computed solution.
We remark that actual calls to \texttt{driver} may produce different results,
because different realizations of the noise vector propagate differently in the
computation.

\begin{figure}[ht]
\begin{center}
\includegraphics[width=0.65\textwidth]{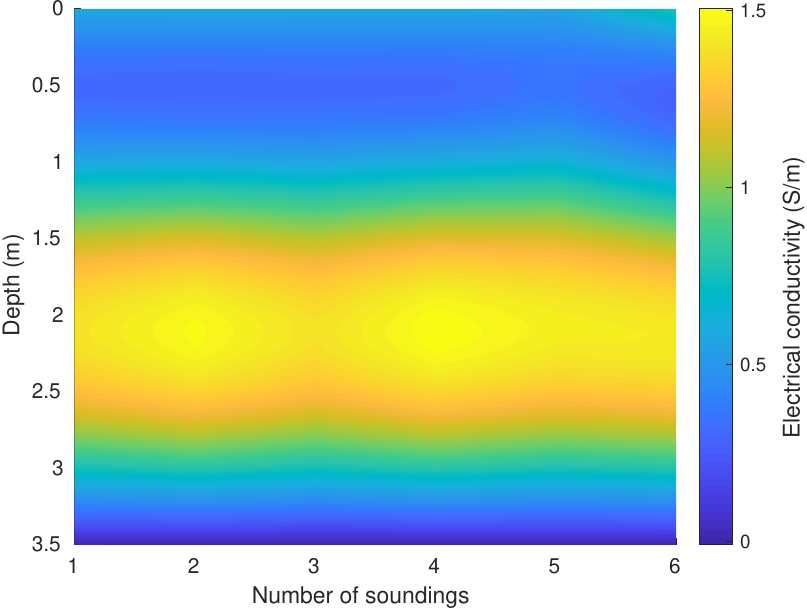}
\end{center}
\caption{Reconstruction of the electrical conductivity from experimental data,
as performed by the \texttt{driver2D} test script.}
\label{Fig2}
\end{figure}

A second test program \texttt{driver2D} is provided, which performs the
processing of an experimental data set composed by 6 measurements.
An additional figure displays the computed solutions in a 2D pseudo-color plot;
see Figure~\ref{Fig2}.

Table~\ref{tab:opts} reports the options that can be set to tune the
functionality of the package. All the options have a default value; see the
\texttt{driver} program for an example of their use, as well as the help pages
of the computational routines.

Further information on the use of the package, as well as on the data structure
adopted to read experimental data sets, is provided in the User Manual,
included in the package main directory.

\begin{table}[ht!]
\centering
\footnotesize
\caption{Options for the tuning of the computational routines in the package.}
\begin{tabular}{p{2.4cm}p{8.5cm}}
\hline
\multicolumn{2}{c}{\textbf{Problem and Data Definition}} \\
\hline
\texttt{func} & function to be minimized, defaults to \texttt{aconduct} \\
\texttt{orientation} & orientation of device (vertical, horizontal, or both) \\
\texttt{yrows} & data rows being processed (useful to exclude some data) \\
\hline
\multicolumn{2}{c}{\textbf{Iterative Method Tuning}} \\
\hline
\texttt{sigmainit} & starting vector for conductivity sigma
	(\texttt{emsolvenlsig}) \\
\texttt{mu} & estimated value of mu in each layer (\texttt{emsolvenlsig}) \\
\texttt{muinit} & starting vector for permeability mu (\texttt{emsolvenlmu}) \\
\texttt{sigma} & estimated value of sigma in each layer 
	(\texttt{emsolvenlmu}) \\
\texttt{tau} & stop tolerance \\
\texttt{nmax} & maximum number of iterations \\
\texttt{damped} & set to 1 for damped method, to 0 for undamped \\
\texttt{dampos} & set to 1 for positive solution, to 0 for unconstrained \\
\texttt{metjac} & method to compute/approximate the Jacobian \\
\texttt{kbroyden} & interval for Broyden's Jacobian updates (see
	\cite{dfr14}) \\
\hline
\multicolumn{2}{c}{\textbf{Regularization Matrix and Parameter Estimation}} \\
\hline
\texttt{Lind} & index of regularization matrix ($k$th derivative) \\
\texttt{metk} & methods for choosing the regularization parameter \\
\texttt{ds} & standard deviation of noise for discrepancy principle \\
\texttt{taudiscr} & $\tau$ factor for discrepancy principle \\
\hline
\multicolumn{2}{c}{\textbf{Displaying Intermediate Results}} \\
\hline
\texttt{show} & if $>0$ show information on iterations \\
\texttt{showpar} & if $>0$ show information on identification of reg.~parameter \\
\texttt{rowcol} & size of grid in the iteration graphical window \\
\texttt{xlim} & x-limit for the graphs in the iteration window \\
\texttt{ylim} & y-limit for the graphs in the iteration window \\
\texttt{truesol} & exact solution, if available (to compute errors) \\
\hline
\end{tabular}
\label{tab:opts}
\end{table}

\begin{figure}[ht]
\begin{center}
\includegraphics[width=0.8\textwidth]{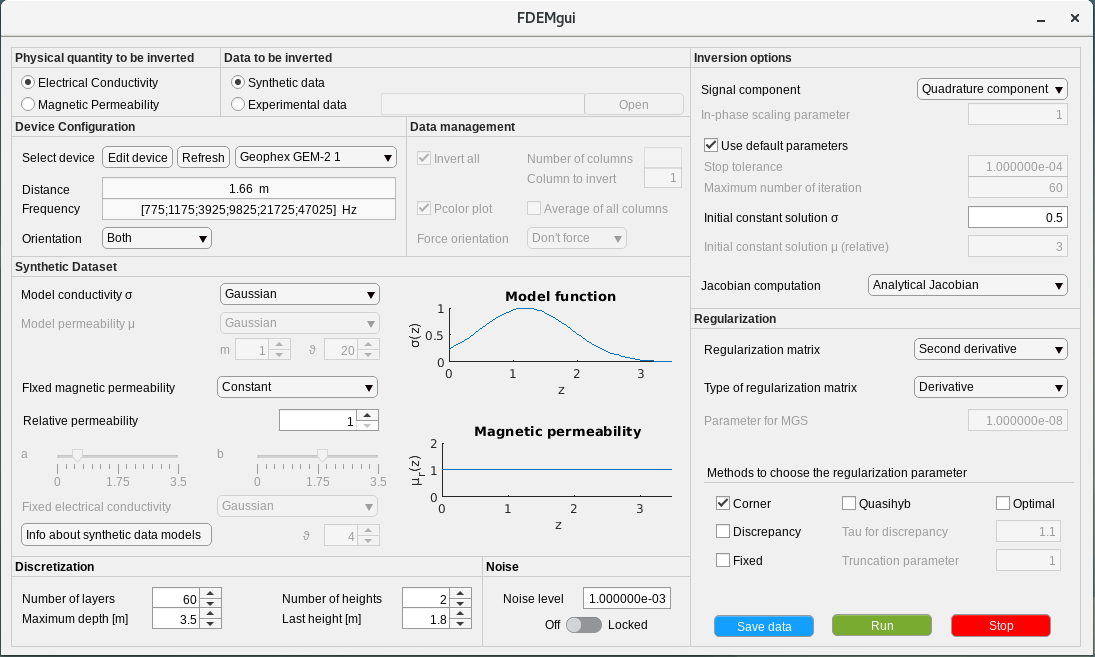}
\end{center}
\caption{The FDEMgui graphical user interface.}
\label{gui}
\end{figure}

While for processing large experimental data sets it is probably preferable to
use a script resembling the \texttt{driver} program,
the toolbox includes a MATLAB graphical user interface (GUI)
in order to simplify experimenting with the functions of the package.
It is started by issuing the command \texttt{FDEMgui} in the main MATLAB
window; see Figure~\ref{gui}.
It is composed of a set of input panels that lets the user generate or load
data sets, choose different approaches and parameters for the inversion
algorithm and, finally, visualize the computed results. 
The panels, which are extensively described in the User Manual, are the
following:
\begin{itemize}
\item Physical quantity to be inverted,
\item Data to be inverted,
\item Device configuration,
\item Data management,
\item Synthetic dataset,
\item Discretization,
\item Noise,
\item Inversion options,
\item Regularization.
\end{itemize}
The GUI is completed by three buttons that allow the user to start the
computation, interrupt it, in case something goes wrong, and save the computed
solution to a \texttt{.mat} file.

\section{Conclusions}
\label{sec:4}

In this paper we presented a new MATLAB toolbox for reconstructing the
electromagnetic features of the subsoil starting from FDEM data sets, by a
regularized, damped, Gauss--Newton method. The package has been developed and
extensively tested in various papers recently published.
The toolbox includes both test programs, that demonstrates its use,
and a graphical user interface, which simplifies setting the parameters
that control the package functions.

\section*{Acknowledgements}
Research partially supported by 
the Fondazione di Sardegna 2017 research project ``Algorithms for Approximation
with Applications [Acube]'',
the INdAM-GNCS research project ``Metodi numerici per problemi mal posti'',
the INdAM-GNCS research project ``Discretizzazione di misure, approssimazione 
di operatori integrali ed applicazioni'', and the Regione Autonoma della Sardegna research project ``Algorithms and Models for Imaging Science [AMIS]''
(RASSR57257, intervento finanziato con risorse FSC 2014-2020 - Patto per lo
Sviluppo della Regione Sardegna). 
CF gratefully acknowledges Regione Autonoma della Sardegna for the financial
support provided under the Operational Programme P.O.R. Sardegna F.S.E.
(European Social Fund 2014-2020 - Axis III Education and Formation, Objective
10.5, Line of Activity 10.5.12).

\bibliographystyle{siam}

\end{document}